\documentclass[journal,transmag]{IEEEtran}
\usepackage{ifpdf}

\usepackage{cite}

\ifCLASSINFOpdf
  \usepackage[pdftex]{graphicx}
\else
   \usepackage[dvips]{graphicx}
\fi

\usepackage{amsmath}
\usepackage{amssymb}
\usepackage{mathtools}

\usepackage{algorithmic}

\usepackage{array}

\ifCLASSOPTIONcompsoc
 \usepackage[caption=false,font=normalsize,labelfont=sf,textfont=sf]{subfig}
\else
 \usepackage[caption=false,font=footnotesize]{subfig}
\fi

\usepackage{stfloats}

\ifCLASSOPTIONcaptionsoff
 \usepackage[nomarkers]{endfloat}
\let\MYoriglatexcaption\caption
\renewcommand{\caption}[2][\relax]{\MYoriglatexcaption[#2]{#2}}
\fi

\usepackage{url}
\usepackage{hyperref}
\usepackage{eso-pic}
\AddToShipoutPicture*{\footnotesize\sffamily\raisebox{0.75cm}{\hspace{1.35cm}\fbox{\parbox{\textwidth}{\copyright\,2023 IEEE.  Personal use of this material is permitted.  Permission from IEEE must be obtained for all other uses, in any current or future media, including reprinting/republishing this material for advertising or promotional purposes, creating new collective works, for resale or redistribution to servers or lists, or reuse of any copyrighted component of this work in other works.}}}}

\usepackage{blindtext}
\usepackage{tikz-cd}
\usepackage{adjustbox}
\usepackage{pgfplots}
\usetikzlibrary{decorations.markings,arrows}
\pgfplotsset{compat=newest}
\usepackage{color}

\begin{document}


\title{Finite Element Modeling of Power Cables \\ using Coordinate Transformations}


\author{\IEEEauthorblockN{Albert Piwonski\IEEEauthorrefmark{1},
Julien Dular\IEEEauthorrefmark{2},
Rodrigo Silva Rezende\IEEEauthorrefmark{1}, 
Rolf Schuhmann\IEEEauthorrefmark{1}}
\IEEEauthorblockA{\IEEEauthorrefmark{1}Theoretische Elektrotechnik,
Technische Universität Berlin, Berlin, Germany}
\IEEEauthorblockA{\IEEEauthorrefmark{2} TE-MPE-PE, CERN, Geneva, Switzerland}%
}


\markboth{COMPUGMAG 2023: OC1: Static and quasi-static fields / Wave propagation. Conference Paper Number: 258.}%
{Shell \MakeLowercase{\textit{et al.}}: Bare Demo of IEEEtran.cls for IEEE Transactions on Magnetics Journals}


\IEEEtitleabstractindextext{%
\begin{abstract}

Power cables have complex geometries in order to reduce their ac resistance. Although there are many different cable designs, most have in common that their inner conductors' cross-section is divided into several electrically insulated conductors, which are twisted over the cable’s length (helicoidal symmetry). In previous works, we presented how to exploit this symmetry by means of dimensional reduction within the $\mathbf{H}-\varphi$ formulation of the eddy current problem. Here, the dimensional reduction is based on a coordinate transformation from the Cartesian coordinate system to a helicoidal coordinate system. This contribution focuses on how this approach can be incorporated into the magnetic vector potential based $\mathbf{A}-v$ formulation. 
\end{abstract}

}

\maketitle
\IEEEdisplaynontitleabstractindextext
\IEEEpeerreviewmaketitle


\section{Introduction}
\label{sec:introduction}

\begin{figure}[t]
    \centering
    \includegraphics[width=0.4\textwidth]{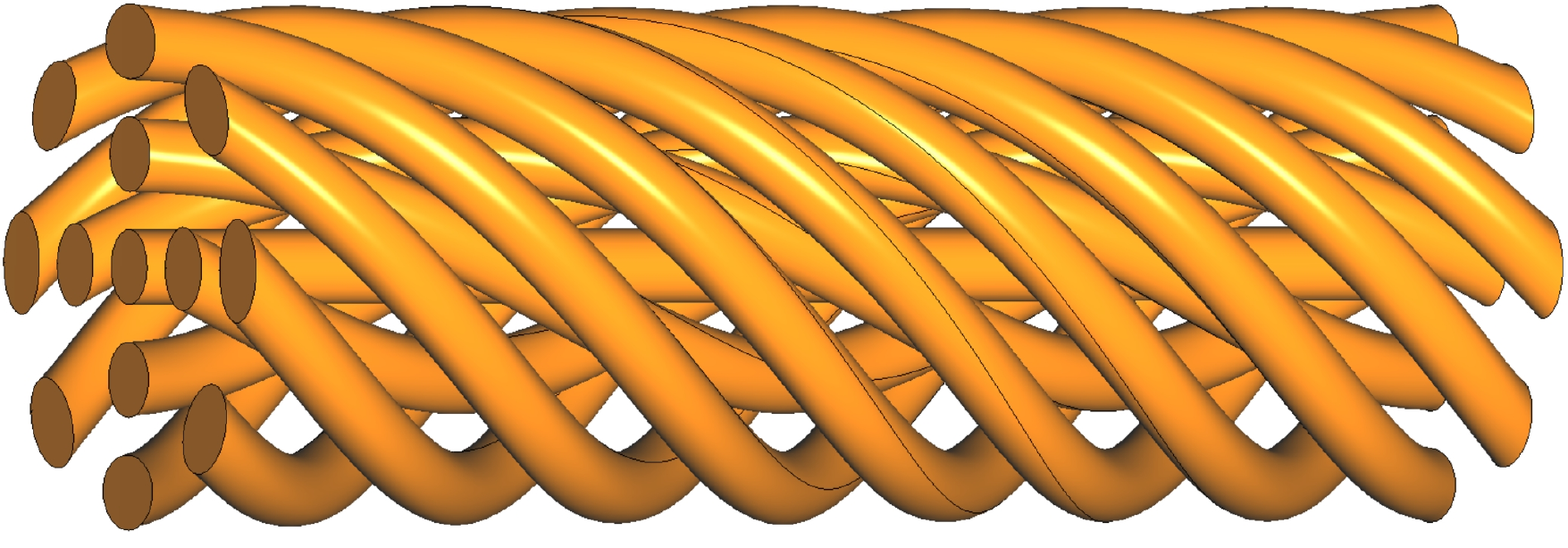}
    \caption{Generic model of a cable's inner conductor: Modeled using \cite{a1}.}
    \label{fig:cst_cable}
\end{figure}

\IEEEPARstart{P}{ower} cables are essential elements in the transmission
chain of electric power from generator to consumer. Special inner conductor designs have been developed for ac operation to minimize the undesirable current displacement caused by the skin and proximity effect~\cite{a2}. Although there are many different cable designs, most have in common that their inner conductors' cross-section is devided into several conductors, which are twisted and electrically insulated from each other (see Fig.~\ref{fig:cst_cable}). 

Investigating ac losses in power cables is an important but challenging task: Solving numerically eddy current boundary value problems (BVP) in 3-D, which model the cable’s electromagnetic
behaviour in the magnetoquasistatic limit, leads to tremendous
computational efforts, due to arising multiscale problems (e.g., thin insulation layers).

However, if one models the power cable as a symmetric BVP, computational costs can be
scaled down significantly by means of dimensional reduction. Due to the conductors' twist, however, no conventional translational symmetry is valid, as recently shown in a formal framework using Lie derivatives~\cite{a222}. In particular though, when choosing proper boundary conditions, an eddy current BVP
posed on a domain as in Fig.~\ref{fig:cst_cable} has a so-called helicoidal symmetry. In contrast to applying periodic boundary conditions, here the
model can be solved in 2-D. This special symmetry was exploited before to calculate hysteresis losses in twisted superconductors~\cite{a4, a22}, but also for solving full Maxwell's equations in twisted waveguides~\cite{a3}. Further, in previous works, we presented how to incorporate this symmetry property within the $\mathbf{H}-\varphi$ finite element formulation applied to a power cable problem~\cite{a23}. 

This contribution focuses on an integration into the magnetic vector potential based $\mathbf{A}-v$ formulation: In Section~\ref{sec:coordinate_transformations}, we define the coordinate transformation and its inverse that allow a bidirectional transition between Cartesian and helicoidal coords. Subsequently, Section~\ref{sec:A_v_formulation_in_different_coordinates} presents how a \mbox{3-D} $\mathbf{A}-v$ formulation in Cartesian coords can be reformulated equivalently into a 2-D formulation in helicoidal coords. Further, implementation details of the 2-D model are given. In Section~\ref{sec:results}, we compare the results of this contribution and of our previous work~\cite{a23} with a 3-D reference model. 

\section{Coordinate transformations}
\label{sec:coordinate_transformations}

In computational electromagnetism, symmetries occur quite differently: E.g., when computing electrical machines, often a translational symmetry is assumed. Further, azimuthal symmetries appear frequently when modeling cylindrical waveguides. The common key idea is the use of a coordinate system in which the geometries appear the same in one direction. Then, partial derivatives w.r.t. to this particular direction of electromagnetic field quantities are either vanishing or constant, leading to drastically reduced computational efforts. Therefore, for describing twisted conductors, a helicoidal coordinate sytem is the most suitable.

In the following, we denote points represented in the Cartesian coordinate system $(x,\,y,\,z)$ as $\mathbf{p}_{xyz}\coloneqq [x,\,y,\,z]^\top$, whereas points represented in the helicoidal coordinate system $(u,\,v,\,w)$ are denoted as $\mathbf{p}_{uvw}\coloneqq [u,\,v,\,w]^\top$. The birectional change of coordinates is achieved by the map $\boldsymbol\phi: \Omega_{xyz} \rightarrow \Omega_{uvw}$ and its inverse $\boldsymbol\phi^{-1}:\Omega_{uvw} \rightarrow \Omega_{xyz}$, where $\Omega_{xyz},\, \Omega_{uvw} \subset \mathbb{R}^3$: 

\begin{align}
    \boldsymbol\phi(\mathbf{p}_{xyz}) &= \mathbf{p}_{uvw} = 
    \begin{bmatrix}
        + x \cos(z\alpha/\beta) + y \sin(z\alpha/\beta) \\
        - x \sin(z\alpha/\beta) + y \cos(z\alpha/\beta) \\
        +z
    \end{bmatrix}, \label{eq:xyz_2_uvw}
    \end{align}
    \begin{align}
    \boldsymbol\phi^{-1}(\mathbf{p}_{uvw}) &= \mathbf{p}_{xyz} =
    \begin{bmatrix}
        +u\cos(w\alpha/\beta) -v\sin(w\alpha/\beta) \\
        +u\sin(w\alpha/\beta) +v\cos(w\alpha/\beta) \\
        +w
    \end{bmatrix}\label{eq:uvw_2_xyz}.
\end{align}

\begin{figure}[h]
    \centering
    \includegraphics[width=0.48\textwidth]{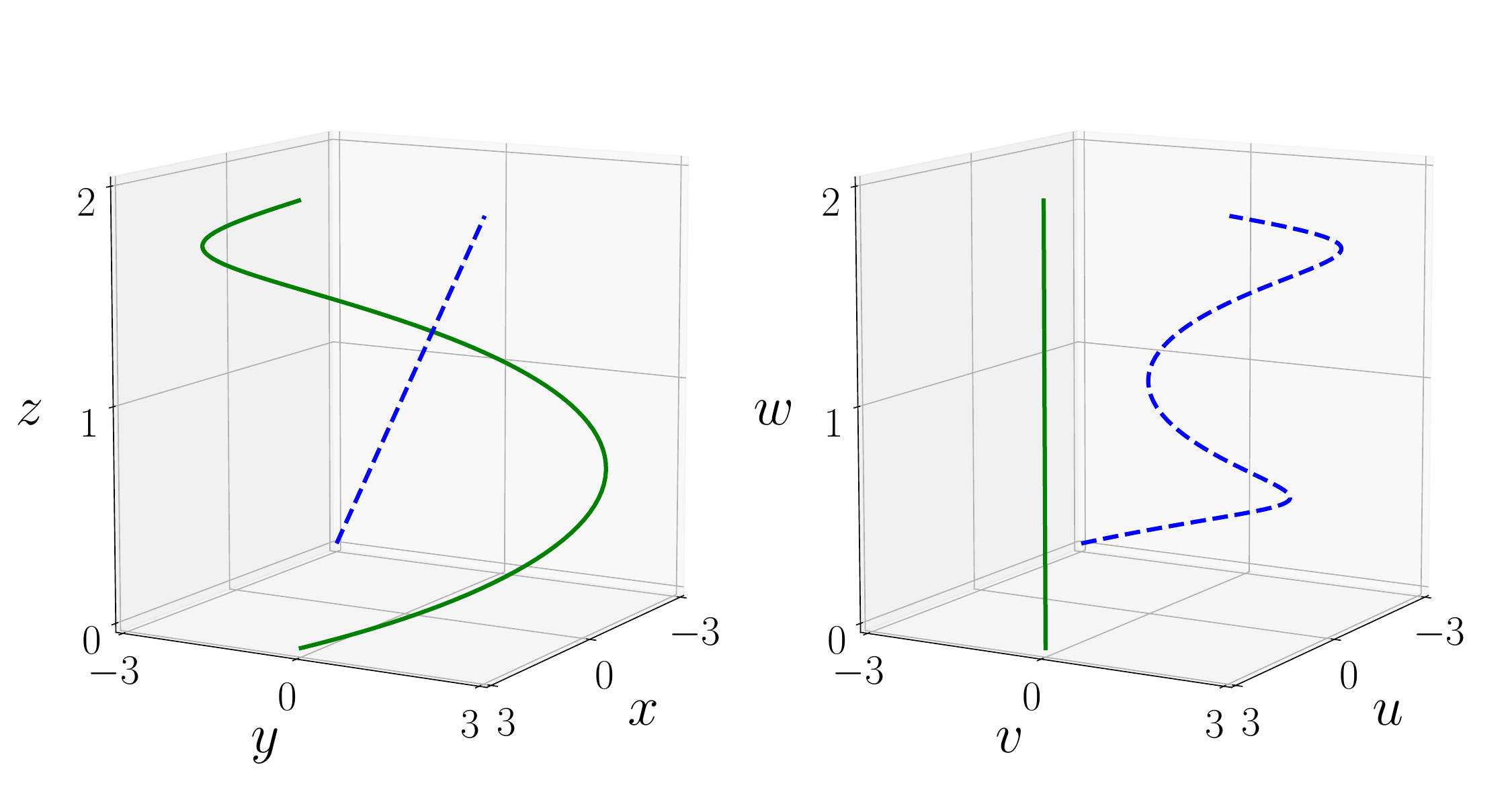}
    \caption{Geometric objects represented in Cartesian (left, arbitrary units) \& helicoidal coordinates (right, same arbitrary units): Helical objects appear straight ($w$-invariant), which is in general not true for straight lines (e.g., see the dashed line between $[-3,\,-3,\,0]^\top$ and $[3,\,3,\,2]^\top$). Further, note that the $(u,\,v,\,w)$ coordinate system must be understood here as non-directional since it is non-orthogonal ($w$-axis is not orthogonal to the $uv$-plane).}
    \label{fig:understanding_transformation}
\end{figure}

\noindent Here, parameters $\alpha$, $\beta$ are related to the number of turns and to the total longitudinal length of the helical object of interest, i.e., for different geometries $\boldsymbol\phi$ resp. $\boldsymbol\phi^{-1}$ are defined differently as well. The effect of the global coordinate transformation is demonstrated in Fig. \ref{fig:understanding_transformation}.

\vspace{-0.45cm}

\begin{figure}[h]
\centering\includegraphics[width=0.27\textwidth]{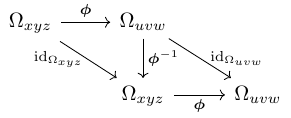}
    \caption{Coordinate transformation $\boldsymbol{\phi}$, its inverse $\boldsymbol{\phi}^{-1}$ and identities on $\Omega_{xyz}$ resp. $\Omega_{uvw}$ summarized as commutative diagram, i.e., $\boldsymbol\phi^{-1} \circ \boldsymbol\phi \equiv \mathrm{id}_{\Omega_{xyz}}$ and $\boldsymbol\phi \circ \boldsymbol\phi^{-1} \equiv \mathrm{id}_{\Omega_{uvw}}$.}
    \label{fig:commutative_diagram}
\end{figure}

It is important to note, that the maps \eqref{eq:xyz_2_uvw} \& \eqref{eq:uvw_2_xyz} allow a unique bidirectional transition of points between both coordinate systems, s.t. the diagram in Fig.~\ref{fig:commutative_diagram} commutes. This property ensures that also the electromagnetic field quantities have a unique representation in the other coordinate system, which will be needed in the following.      

\section{\texorpdfstring{$\mathbf{A}-v$}{TEXT} formulation in different coordinates}
\label{sec:A_v_formulation_in_different_coordinates}

In previous works, we presented how to exploit helicoidal symmetries, in the context of eddy current problems, by means of dimensional reduction within the \mbox{$\mathbf{H}-\varphi$} finite element formulation~\cite{a23}.  In this paper, on the other hand, we present how to incorporate this approach into the more widespread \mbox{$\mathbf{A}-v$} formulation. 

From a topological point of view, the latter formulation is simpler, since typically no topological pre-processing of the computational domain~$\Omega$ is needed. In a nutshell, this can be explained by the fact that the magnetic flux density $\mathbf{B}$ in $\Omega$ is closed ($\mathrm{div}\,\mathbf{B}=0$) and exact ($\exists\,\mathbf{A}$, s.t. $\mathbf{B} = \mathbf{curl}\,\mathbf{A}$)~\cite{lindell}. In contrast, within the \mbox{$\mathbf{H}-\varphi$} formulation, we consider the magnetic field $\mathbf{H}$ in the insulating domain $\Omega_i$, which is also closed ($\mathbf{curl}\,\mathbf{H} = \mathbf{0}$) as there are no currents, but not exact ($\nexists\,\varphi$, s.t. $\mathbf{grad}\,\varphi = \mathbf{H}$). This means in particular, if one represents the magnetic field purely as the gradient of a scalar potential $\varphi$, Ampere's law would be violated in the presence of conductors, since then no net circulations of $\mathbf{H}$ along closed curves could be generated.


\subsection{Weak formulation in Cartesian coordinates}
\label{subsec:Weak_formulation_in_Cartesian_coordinates}

In the following, we present the weak \mbox{$\mathbf{A}-v$} formulation in a concise manner. More details can be found in~\cite{a22}. In the Cartesian domain $\Omega_{xyz}$ with boundary $\partial\Omega_{xyz}$, conducting subdomain $\Omega_{xyz,c}$ ($N$ connected components), and in frequency domain with frequency $\omega/2\pi$, the \mbox{$\mathbf{A}-v$} formulation consists of two sets of equations:

{\small
\begin{align}
    &\left(\nu_{xyz}\mathbf{curl}\, \mathbf{A}_{xyz},\,\mathbf{curl}\,\mathbf{A}_{xyz}'  \right)_{\Omega_{xyz}} \nonumber \\ + &\left(j\omega\sigma_{xyz}\,\mathbf{A}_{xyz},\,\mathbf{A}_{xyz}'  \right)_{\Omega_{xyz,c}} -  
    \langle\widetilde{\mathbf{H}}_{xyz} \times \mathbf{n},\, \mathbf{A}_{xyz}'\rangle_{\partial\Omega_{xyz}} \nonumber \\
       + &\left(\sigma_{xyz}\,\mathbf{grad}\,{v}_{xyz},\,\mathbf{A}_{xyz}'  \right)_{\Omega_{xyz,c}}  = 0, \label{eq:a_v_weak_part_1} \\ \nonumber \\
   &\left(j\omega\sigma_{xyz}\,\mathbf{A}_{xyz},\,\mathbf{grad}\,v_{xyz}'  \right)_{\Omega_{xyz,c}} \nonumber \\ + &\left(\sigma_{xyz}\,\mathbf{grad}\,v_{xyz},\,\mathbf{grad}\,v_{xyz}'  \right)_{\Omega_{xyz,c}}= \sum_{i=1}^N I_i V_i',\label{eq:a_v_weak_part_2}
\end{align}
}

\noindent with the magnetic vector potential $\mathbf{A}_{xyz}$ and electric scalar potential $v_{xyz}$ (with test functions $\mathbf{A}_{xyz}'$, resp. $\mathbf{grad}\,{v}_{xyz}'$), magnetic reluctivity~$\nu_{xyz}$ and electric conductivity $\sigma_{xyz}$. The term $\widetilde{\mathbf{H}}_{xyz} \times \mathbf{n}$ denotes the fixed tangential magnetic field at the boundary of the domain. Further, the notations $(\mathbf{f}, \, \mathbf{g})_{\Omega_{xyz}}$ and $\langle\mathbf{f}, \, \mathbf{g}\rangle_{\partial\Omega_{xyz}}$ denote $\int_{\Omega_{xyz}} \mathbf{f} \cdot \mathbf{g}\,\mathrm{d}V_{xyz}$ and $\int_{\partial\Omega_{xyz}} \mathbf{f} \cdot \mathbf{g}\,\mathrm{d}S_{xyz}$. 

The electric scalar potential $v_{xyz}$ is only defined in the conducting domain and is further used to associate global quantities to the $i$-th conductor, namely, its total current $I_i$ and its voltage drop $V_i$~\cite{a22}. In $\Omega_{xyz,c}$, the magnetic vector potential is directly linked to the physical electric field $\mathbf{E}_{xyz} = -j\omega\mathbf{A}_{xyz} - \mathbf{grad}\,v_{xyz}$, s.t. here no gauging is needed (modified vector potential). However, in the insulating domain~$\Omega_{xyz,i}$, gauging is necessary to find a unique solution.


\subsection{Weak formulation in helicoidal coordinates}

The key idea for the dimensional reduction is the following: Instead of computing the integrals in~\eqref{eq:a_v_weak_part_1} \& \eqref{eq:a_v_weak_part_2} on the Cartesian domain~$\Omega_{xyz}$, we transform them on the domain $\Omega_{uvw}$, i.e., we shift their computation to a $w$-invariant space (see Fig.~\ref{fig:understanding_transformation}). 

This transformation requires a reexpression of the involved electromagnetic field quantities. The one-forms $\mathbf{A}_{xyz}$, $\mathbf{grad}\,v_{xyz}$ and $\widetilde{\mathbf{H}}_{xyz}$ and the two-form $\mathbf{curl}\, \mathbf{A}_{xyz}$ (magnetic flux density $\mathbf{B}_{xyz}$), transform as follows~\cite{a22, lindell}:
\begin{align}
    \mathbf{A}_{xyz}(\mathbf{p}_{xyz}) &= J_{\boldsymbol\phi^{-1}}^{-\top}\ \mathbf{A}_{uvw}(\boldsymbol\phi(\mathbf{p}_{xyz})),\label{eq:A_uvw_2_A_xyz}\\
    \ \mathbf{curl}\,\mathbf{A}_{xyz}(\mathbf{p}_{xyz}) &= \frac{J_{\boldsymbol\phi^{-1}}}{\mathrm{det}(J_{\boldsymbol\phi^{-1}})}\ \mathbf{curl}\,  \mathbf{A}_{uvw}(\boldsymbol\phi(\mathbf{p}_{xyz})),\label{eq:B_uvw_2_B_xyz}
\end{align}
where $J_{\boldsymbol\phi^{-1}}$ denotes the Jacobian of $\boldsymbol\phi^{-1}$ evaluated at point $\mathbf{p}_{uvw}$. Here, it is important to mention, that the $\mathbf{curl}$ operator on the right hand side of eq.~\eqref{eq:B_uvw_2_B_xyz} is not the actual differential operator in the helicoidal coordinate system. Rather, this operator is to be understood as blindly applied (applied as in Cartesian coords with a relabeling of variables and components) -- its results are then directly mapped back into the Cartesian coordinate system as the equality sign in eq.~\eqref{eq:B_uvw_2_B_xyz} holds componentwise for vector fields expressed in Cartesian coordinates. Also, the qualitatively same statement applies, when inserting $\mathbf{grad}\,v_{xyz}$ into formula~\eqref{eq:A_uvw_2_A_xyz}. Changing variables also introduces, rather conceptually than computationally, a factor $\mathrm{det}(J_{\boldsymbol\phi^{-1}})= 1$ in the integrals in eq.~\eqref{eq:a_v_weak_part_1} \& \eqref{eq:a_v_weak_part_2}. \ \\

In terms of the $(u,\,v,\,w)$ coordinates, introducing trial function spaces $A(\Omega_{uvw})$ \& $V(\Omega_{uvw})$ and test function spaces $A_0(\Omega_{uvw})$ \& $V_0(\Omega_{uvw})$ which will be defined in Section~\ref{subsec:Space_discretization_and_implementation_details}, we can reformulate the weak formulation as follows: Seek $\mathbf{A}_{uvw} \in A(\Omega_{uvw})$ and  $\mathbf{grad}\, v_{uvw} \in V(\Omega_{uvw,c})$, s.t. $\forall \,\mathbf{A}_{uvw}'\in A_0(\Omega_{uvw})$ and $\forall \,\mathbf{grad}\, v_{uvw}'\in V_0(\Omega_{uvw,c}) $:

{\small
\begin{align}
    &\left(\boldsymbol\nu_{uvw}\mathbf{curl}\, \mathbf{A}_{uvw},\,\mathbf{curl}\,\mathbf{A}_{uvw}'  \right)_{\Omega_{uvw}} \nonumber \\ + &\left(j\omega\boldsymbol\sigma_{uvw}\,\mathbf{A}_{uvw},\,\mathbf{A}_{uvw}'  \right)_{\Omega_{uvw,c}} -  
    \langle\widetilde{\mathbf{H}}_{xyz} \times \mathbf{n},\, \mathbf{A}_{uvw}'\rangle_{\partial\Omega_{uvw}} \nonumber \\
       + &\left(\boldsymbol\sigma_{uvw}\,\mathbf{grad}\,{v}_{uvw},\,\mathbf{A}_{uvw}'  \right)_{\Omega_{uvw,c}}  = 0, \label{eq:a_v_uvw_weak_part_1} \\ \nonumber \\
   &\left(j\omega\boldsymbol\sigma_{uvw}\,\mathbf{A}_{uvw},\,\mathbf{grad}\,v_{uvw}'  \right)_{\Omega_{uvw,c}} \nonumber \\ + &\left(\boldsymbol\sigma_{uvw}\,\mathbf{grad}\,v_{uvw},\,\mathbf{grad}\,v_{uvw}'  \right)_{\Omega_{uvw,c}}= \sum_{i=1}^N I_i V_i',\label{eq:a_v_uvw_weak_part_2}
\end{align}
}

\noindent where the effect of the change of variables is fully contained in two anisotropic, $w$-invariant material parameters, written as tensors:

\begin{align}
    \boldsymbol\sigma_{uvw}(\mathbf{p}_{uvw}) &= \sigma_{xyz}(\boldsymbol\phi^{-1}(\mathbf{p}_{uvw})) J_{\boldsymbol\phi^{-1}}^{-1}J_{\boldsymbol\phi^{-1}}^{-\top}\mathrm{det}(J_{\boldsymbol\phi^{-1}}), \label{mu_xyz_2_uvw} \\
    \boldsymbol\nu_{uvw}(\mathbf{p}_{uvw}) &= \nu_{xyz}(\boldsymbol\phi^{-1}(\mathbf{p}_{uvw})) J_{\boldsymbol\phi^{-1}}^{\top}J_{\boldsymbol\phi^{-1}}/\mathrm{det}(J_{\boldsymbol\phi^{-1}}). \label{rho_xyz_2_uvw} 
\end{align}

\noindent Now, since neither the integrands nor the computational domain~$\Omega_{uvw}$ depend on the $w$-coordinate, one may solve the problem on any $uv$-plane for a fixed value of $w$. For simplicity, we choose $w = 0$, because then $x=u$ and $y=v$, see eq.~\eqref{eq:uvw_2_xyz}.


\subsection{Space discretization and implementation details}
\label{subsec:Space_discretization_and_implementation_details}

We solve the systems of equations~\eqref{eq:a_v_uvw_weak_part_1} \& \eqref{eq:a_v_uvw_weak_part_2} by using the open-source
finite-element software GetDP~\cite{d1}, which allows for flexible
function space definitions, whereas the meshing process is
performed by Gmsh controlled via the Julia API~\cite{d2, c4}. In previous works, we derived the \mbox{2-D} computational domain (symmetry cell) analytically using the theory of envelopes~\cite{a23}. Now, the cross-section is computed as the intersection of 3-D helicoidal conductors with a plane using built-in CAD routines. This will allow us to investigate elliptically shaped conductor cross-sections in future work, which can be created under mechanical stress in the cable manufacturing process~\cite{a2}.

Although we solve the problem on a 2-D domain, we still assume, that the magnetic vector potential~$\mathbf{A}_{uvw}$ has full three components, which are interpolated separately: The in-plane components ($A_u$ \& $A_v$) are discretized using 2-D Whitney edge functions. In the domain~$\Omega_{uvw,i}\setminus\partial\Omega_{uvw,c}$, though, only degrees of freedom (DOF) associated with the cotree of the mesh are non-zero, i.e., a tree-cotree gauge is applied to achieve a unique solution~\cite{c5}. This defines the function space $A_{uv}(\Omega_{uvw})$. 

Further, the component $A_w$ is spanned by nodal functions and is essentially fixed to zero on $\partial\Omega_{uvw}$ as we model the cable's shield as a perfect electrically conductive cylinder. Since we do not introduce an explicit gauge for $A_w$, this component is implicitly Coulomb-gauged. This constitutes the function space $A_{w}(\Omega_{uvw})$, s.t. in total $A(\Omega_{uvw}) = A_{uv}(\Omega_{uvw}) \oplus A_{w}(\Omega_{uvw})$. The function $\mathbf{grad}\,v_{uvw} \in V(\Omega_{uvw})$ consists of $N$ $w$-directed constants (one constant per connected component of $\Omega_{uvw,c}$).


\section{Results}
\label{sec:results}

\begin{figure}[ht]
    \centering
    \includegraphics[width=0.4\textwidth]{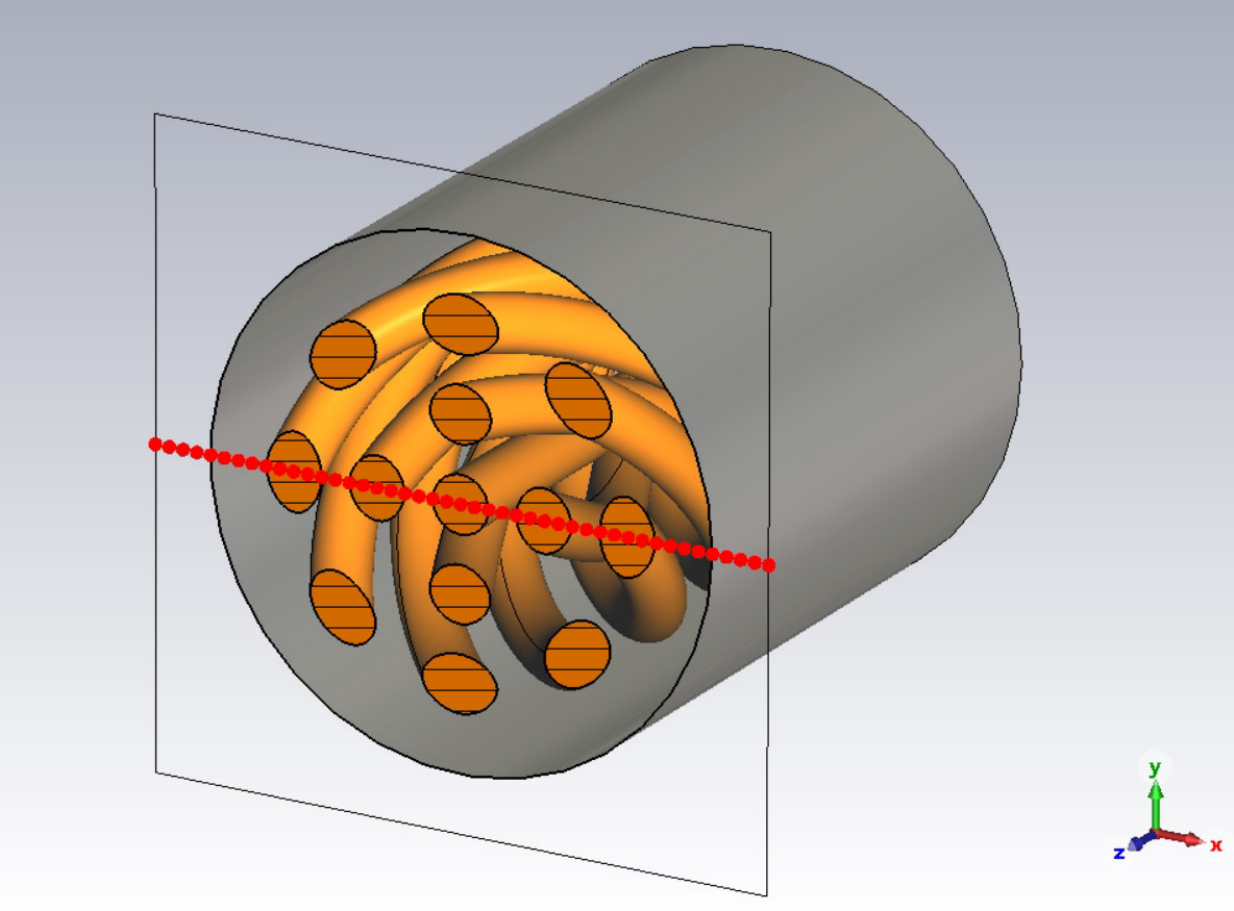}
    \caption{Cross-section view of the 3-D reference model at center $z=0.1\,\mathrm{m}$: Highlighted points along $x$-axis.}
    \label{fig:cst_current_ports}
\end{figure}

The results of the 2-D \mbox{$\mathbf{A}-v$} formulation and of our previously presented 2-D \mbox{$\mathbf{H}-\varphi$} formulation~\cite{a23} are compared with a 3-D reference cable model that is implemented using the commercial software CST Studio Suite~\cite{a1}, referred to as CST. Each of the 13 helicoidal conductors (longitudinal length $0.2\,\text{m}$, cross-section radius $5\,\text{mm}$) carries a total current of amplitude $\sqrt{2}/13 \, \mathrm{A}$ at $f = 50 \, \text{Hz}$. In all models, we considered annealed copper for the conductors' material (conductivity $\approx 58.12 \cdot 10^6 \,\text{S}/\text{m}$) and further assumed a non-magnetic material in the whole domain ($\nu_0 = \mu_0^{-1}$ with $\mu_0 = 4\pi\cdot 10^{-7}\, \text{H}/\text{m}$). Since the current constraint is differently implemented in the 3-D model (current ports located at the cuboid bounding box of the domain), we compare local quantities only at the cable's center (see Fig.~\ref{fig:cst_current_ports}). 

In the following, we compare the 2-D models discretized by $58.64\text{k}$ triangles with a 3-D model discretized by $1.19\text{M}$ tetrahedra: We used first order finite elements resulting in $98\text{k}$ DOF in the 2-D \mbox{$\mathbf{A}-v$} formulation, $57\text{k}$ DOF in the 2-D \mbox{$\mathbf{H}-\varphi$} formulation and $1.21\text{M}$ DOF in the 3-D model. The DOF of the 2-D models differ so much on the same mesh, since only within the \mbox{$\mathbf{H}-\varphi$} formulation we can directly incorporate and exploit the fact that the $w$-component of the magnetic field~$\mathbf{H}_{uvw}$ is a constant in the insulating domain~$\Omega_{uvw,i}$~\cite{a23}. So far, no way has been found to utilize this knowledge also in the \mbox{$\mathbf{A}-v$} formulation. 
The 2-D \mbox{$\mathbf{A}-v$} model outputs the finite-element approximated vector potential $\mathbf{A}_{uvw}$ and the gradient of the electric scalar potential $\mathbf{grad}\, v_{uvw}$. 

Using transformation rules~\eqref{eq:A_uvw_2_A_xyz} \& \eqref{eq:B_uvw_2_B_xyz}, all electromagnetic field quantities in the Cartesian coordinate system can be derived, e.g.:
\begin{align}
    \mathbf{J}_{xyz} &= -\sigma_{xyz} J_{\boldsymbol\phi^{-1}}^{-\top}\,\left(j\omega\mathbf{A}_{uvw} + \mathbf{grad}\,v_{uvw}\right), \label{eq:reconstruction_jxyz}\\
    \mathbf{H}_{xyz} &= \mu_{xyz}^{-1} \frac{J_{\boldsymbol\phi^{-1}}}{\mathrm{det}(J_{\boldsymbol\phi^{-1}})}\ \mathbf{curl}\,  \mathbf{A}_{uvw}\label{eq:reconstruction_hxyz}. 
\end{align}

\begin{figure}[h]
    \centering
    \includegraphics[width=0.5\textwidth]{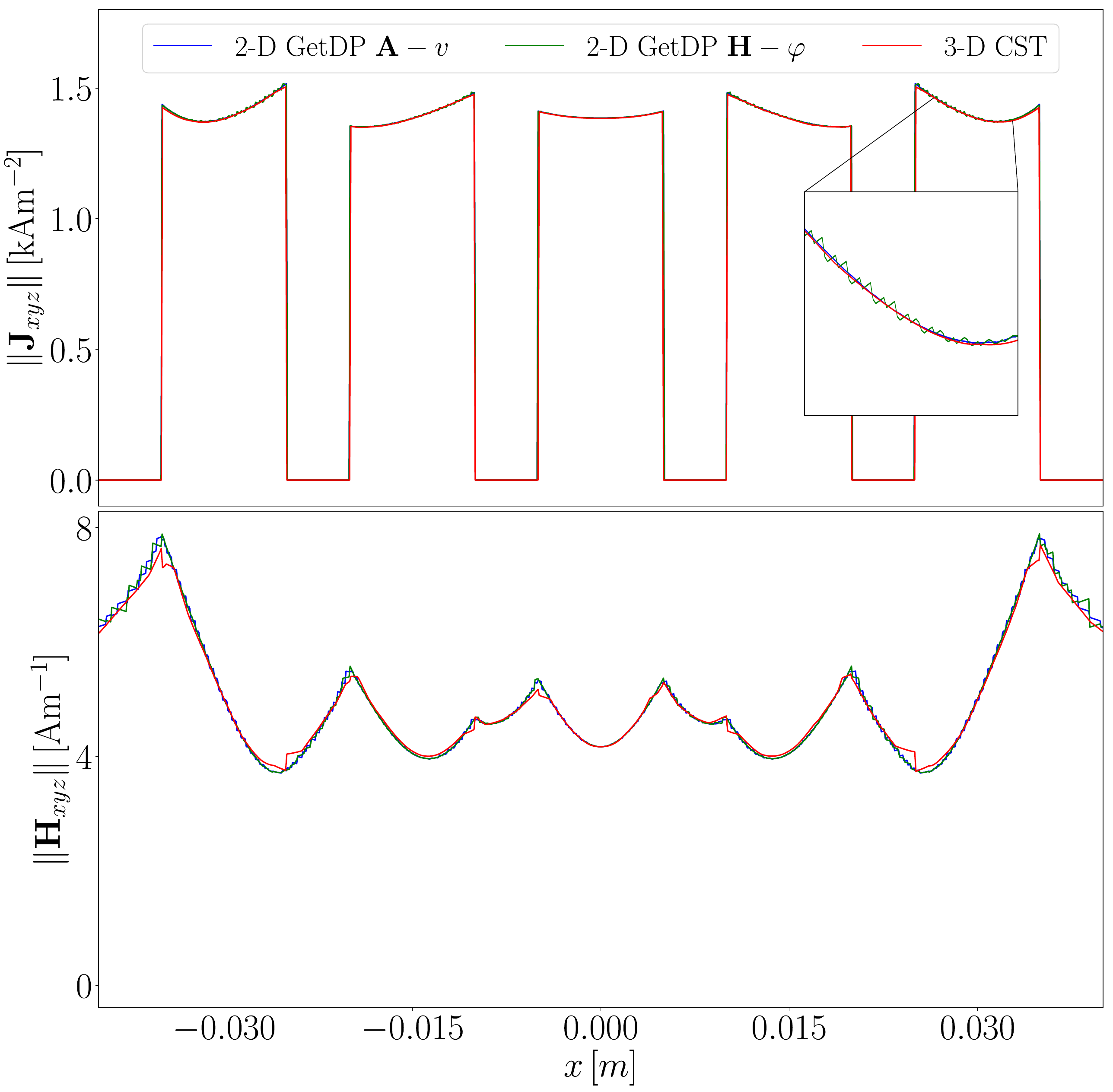}
    \caption{Absolute of $\mathbf{J}_{xyz}$ and $\mathbf{H}_{xyz}$ along $x$-axis.}
    \label{fig:J_xyz_H_xyz}
\end{figure}

\noindent As a local comparison, $\mathbf{J}_{xyz}$ and $\mathbf{H}_{xyz}$ are evaluated along
the $x$-axis. The results depicted in Fig.~\ref{fig:J_xyz_H_xyz} show an accurate agreement between all models. In the \mbox{$\mathbf{A}-v$} formulation, the current density $\mathbf{J}_{xyz}$ is multiplicatively linked to the linear interpolated magnetic vector potential $\mathbf{A}_{uvw}$ (and the piecewise-constant vector field $\mathbf{grad}\, v_{uvw}$), see eq.~\eqref{eq:reconstruction_jxyz}. In contrast, in the \mbox{$\mathbf{H}-\varphi$} formulation, the current density is multiplicatively linked to the $\mathbf{curl}$ operator applied to the linear interpolated magnetic field $\mathbf{H}_{uvw}$~\cite{a23}. Therefore, the current density resulting from the \mbox{$\mathbf{A}-v$} formulation appears comparatively smooth, as shown in the zoom  window in the upper plot of Fig.~\ref{fig:J_xyz_H_xyz}.

Furthermore, the comparison of the ohmic losses, representing
a global quantity, shows a good match: both 2-D models output
a length-related power loss of $21.9\, \mu\text{Wm}^{-1}$, whereas
the 3-D model has a total loss of $4.34\, \mu\text{W}$. Scaling the
length-related losses up to the cable’s longitudinal length results into a
loss of $4.38\, \mu\text{W}$, which deviates 0.9\% from the 3-D result.
We suspect that this discrepancy is mainly due to the different
excitation types.


\newpage
\section{Conclusion \& outlook}
\label{sec:conclusion}

In this work, we have shown how a coordinate transformation can be used for the dimensional reduction of eddy current problems modeling helicoidal symmetric power cables. In particular, we have focused on integrating the symmetry approach into the vector potential based $\mathbf{A}-v$ formulation that can be solved in 2-D (Section~\ref{sec:A_v_formulation_in_different_coordinates}). The results of this 2-D model are in excellent agreement with our previously presented 2-D model~\cite{a23} and with the 3-D reference model (Section~\ref{sec:results}). This reinforces the confidence in the symmetry approach for reducing computational costs when analyzing power cables' electromagnetic behaviour numerically. Future works include the handling of also non-ideal symmetries, see e.g.,~\cite{a22, a222}.





\section*{Acknowledgment}
Special thanks to Christophe Geuzaine (University of Liège) who provided helpful advice on implementations in GetDP.
\ifCLASSOPTIONcaptionsoff
  \newpage
\fi




\begin{thebibliography}{1}



\bibitem{a2} R. Suchantke, ``Alternating Current Loss Measurement of Power Cable Conductors with Large Cross Sections Using Electrical Methods,'' Ph.D. dissertation, Tech. Univ. Berlin, Berlin, Germany, 2018.

\bibitem{a1} \textit{CST Studio Suite}. (2021), Dassault Systèmes, Accessed: Mar. 18, 2022. [Online].
Available:\href{https://www.3ds.com/products-services/simulia/products/cst-studio-suite/}{https://www.3ds.com/products-services/simulia/products/cst-studio-suite/}.

\bibitem{a222}
A. Marjamäki, T. Tarhasaari and P. Rasilo, "Utilizing Helicoidal and Translational Symmetries Together in 2-D Models of Twisted Litz Wire Strand Bundles," \textit{IEEE Transactions on Magnetics}, vol. 59, no. 5, pp. 1-4, May 2023, Art no. 7400504, doi: 10.1109/TMAG.2023.3237767.

\bibitem{a4} A. Stenvall, F. Grilli and M. Lyly, ``Current-Penetration Patterns in Twisted Superconductors in Self-Field,'' \textit{IEEE Transactions on Applied Superconductivity}, 2013, vol. 23, no. 3, pp. 8200105-8200105, Art no. 8200105, doi: 10.1109/TASC.2012.2228733.

\bibitem{a22} J. Dular, ``Standard and Mixed Finite Element Formulations for Systems with Type-II Superconductors,'' Ph.D. dissertation, Univ. of Liège, Liege, Belgium, 2023.

\bibitem{a3} A. Nicolet and F. Zolla, ``Finite element analysis of helicoidal waveguides,'' \textit{IET Science, Measurement \& Technology}, 2007, vol. 1, pp. 67--70,  doi: 10.1049/iet-smt:20060042.

\bibitem{a23} A. Piwonski, J. Dular, R. S. Rezende and R. Schuhmann, "2-D Eddy Current Boundary Value Problems for Power Cables With Helicoidal Symmetry," \textit{IEEE Transactions on Magnetics}, vol. 59, no. 5, pp. 1-4, May 2023, Art no. 6300204, doi: 10.1109/TMAG.2022.3231054.





\bibitem{lindell} I.V. Lindell, ``Differential forms in electromagnetics," John Wiley \& Sons, 2004.

\bibitem{d1} P. Dular, C. Geuzaine, F. Henrotte and W. Legros, ``A general environment for the treatment of discrete problems and its application to the finite element method,'' \textit{IEEE Transactions on Magnetics},  vol.~34, no.~5, pp.~3395--3398, 1998.

\bibitem{d2} C. Geuzaine and J. Remacle, ``Gmsh: a three-dimensional finite element mesh generator with built-in pre- and post-processing facilities,'' \textit{International Journal for Numerical Methods in Engineering}, 2009, vol. 79, no. 11, pp. 1309-1331.

\bibitem{c4} J. Bezanson, A. Edelman, S. Karpinski and V. Shah, ``Julia: A Fresh Approach to Numerical Computing,'' \textit{SIAM Review}, 2017, vol.~59, no.~1, pp.~65-98, doi: 10.1137/141000671.

\bibitem{c5} E. Creusé, P. Dular and S. Nicaise, ``About the gauge conditions arising in finite element magnetostatic problems,'' \textit{Computers \& Mathematics with Applications}, vol.~77, no.~6, pp.~1563-1582, 2019.











\end{thebibliography}
\end{document}